\documentclass[12pt]{article}     \usepackage{amsmath}
\usepackage{amsthm}               \usepackage{amsopn}
\usepackage{psfig}                \usepackage{graphics}
\usepackage{latexsym}             \usepackage{amsfonts}
\usepackage{amssymb}
\swapnumbers
\theoremstyle{plain}
\newtheorem{theorem}{Theorem}[section]

\newtheorem{lemma}[theorem]{Lemma}
\newtheorem{corollary}[theorem]{Corollary}
 
\theoremstyle{definition}

\newtheorem{remark}[theorem]{Remark}
 
\newenvironment{pf}{{\noindent\bf Proof.}}

\newcommand{\nc}[2]{ \newcommand{#1}{#2} }

\nc{\avint}{ {- \hspace{-3.5mm} \int} }  % average integral

\nc{\R}{\rm I \! R}  % Real Numbers
\nc{\N}{\rm I \! N}
\newcommand{\closure}[1]{ \stackrel{\rule{.1 in}{.01 in}}{#1} }

\newcommand{\qclosure}[1]{ \stackrel{\rule{.4 in}{.01 in}}{#1} }

\newcommand{\chisub}[1]{ {\mathbf{\chi}}_{_{#1}} }
 
\newcommand{\newsec}[2]{ \section{#1} \label{sec-#2}  % starts new section,
                         \setcounter{equation}{0}     % resets counters,
                         \setcounter{theorem}{0} }    % makes label

\newcommand{\refeqn}[1]{ (\!\!~\ref{eq:#1}) } % gives references to
\newcommand{\refthm}[1]{ (\!\!~\ref{#1}) }    % equations or theorems
 
\nc{\Holder}{H\"{o}lder\ }

\nc{\ith}{ \ensuremath{\text{i}^{\text{th}}} }
\nc{\jth}{ \ensuremath{\text{j}^{\text{th}}} }
\nc{\kth}{ \ensuremath{\text{k}^{\text{th}}} }
\nc{\curl}{ \nabla \times }
\nc{\Div}{ \nabla \cdot } 

\nc{\Ppl}{ \mathcal{M}^{+} }  \nc{\Pmn}{ \mathcal{M}^{-} }
 
\nc{\smiley}{ $\stackrel{\because}{\smile} \;$ }

\begin{document} 

\numberwithin{equation}{section}

\title{ Convergence of Rothe's Method for\\
       Fully Nonlinear Parabolic Equations }
\author{Ivan Blank and Penelope Smith\\
\normalsize}
 
\maketitle
 
\begin{abstract}
Convergence of Rothe's method for the fully nonlinear parabolic equation
$u_t + F(D^2 u, Du, u, x, t) = 0$ is considered under some continuity
assumptions on $F.$  We show that the Rothe solutions are Lipschitz
in time, and they solve the equation in the viscosity sense.  As an
immediate corollary we get Lipschitz behavior in time of the viscosity
solutions of our equation.
\end{abstract} 

\newsec{Introduction}{AssStruct}
We consider the boundary value problem:
\begin{equation}
\begin{array}{rll}
u_t + F(D_x^2 u, \; D_x u, \; u, \; x, \; t) & \!\!\!= 0 \ \ & \text{in} \ D \\
u & \!\!\!= 0 \ \ & \text{on} \ \partial_{p} D \;. \\
\end{array}
\label{eq:BasicEqn}
\end{equation}
We assume that $u(x,t)$ satisfies Equation\refeqn{BasicEqn}in the viscosity sense,
that $D = \Omega \times [0,T],$ and that $F$ is uniformly elliptic.  We recall that
a function $u$ is a \textit{viscosity subsolution} (\textit{supersolution}) of
Equation\refeqn{BasicEqn}if
   for any $\varphi \in C^{2,1}(\closure{D})$ and any $\epsilon > 0$ which
              satisfies
     \begin{equation}  \hspace{-.1in}
        \begin{array}{ll}
     & \varphi_t + F(D^2 \varphi, \; D \varphi, \; \varphi, \; x, \; t) \geq \epsilon > 0 \\
     (\text{or} &
     \varphi_t + F(D^2 \varphi, \; D \varphi, \; \varphi, \; x, \; t) \leq - \epsilon < 0
     \ \ \text{for the supersolution case})
        \end{array}
     \label{eq:vpissuper}
     \end{equation}
              $u - \varphi$ cannot attain a local maximum (minimum) of $0.$  (In other words,
              $\varphi$ cannot touch $u$ from above (below).)
$u$ is a \textit{solution} if it is both a subsolution and a supersolution.  Uniform
ellipticity means that there are positive numbers $\lambda$ and $\Lambda$ such that
for any positive definite matrix $N$ we have
\begin{equation}
- \Lambda ||N|| \leq F( M + N, \; P, \; v, \; x, \; t)
                 - F( M, \; P, \; v, \; x, \; t) \leq - \lambda ||N|| \;,
\label{eq:UnifEllip}
\end{equation}
or what is equivalent:
\begin{equation}
F( M + N, \; P, \; v, \; x, \; t) \leq
F( M, \; P, \; v, \; x, \; t) + \Lambda ||N^{-}|| - \lambda ||N^{+}||
\label{eq:UnifEllip2}
\end{equation}
for all $M$ and $N.$  So the negative Laplacian would fit our structure conditions,
and not the positive Laplacian.
 
For linear $F$'s Rothe's method is commonly used as a numerical approximation.  Rothe's
method corresponds to doing a backward Euler approximation in Banach space, and is also
known sometimes as the method of lines.  Kikuchi and Ka\v{c}ur have studied the
convergence properties of Rothe's method in a series of recent papers in a variety of
function spaces.  (See \cite{K} and \cite{KK} and the references therein.)  Pluschke has
studied the quasilinear case (see \cite{P}), but the current work appears to be the first
extension of Rothe's method to viscosity solutions in the fully nonlinear setting.
 
We fix our mesh size, $0 < h \leq 1,$ and define $z_{n,h}(x)$ recursively.
We take $z_{0,h}(x) \equiv 0,$ and for $n \geq 0$ we let $z_{n+1,h}$ be the solution
of the elliptic problem:
\begin{equation}
\left.
\begin{array}{rl}
\begin{array}{r}
F(D^2 z_{n+1,h}, \; D z_{n+1,h}, \; z_{n+1,h}, \; x, \; (n+1)h) +
   \frac{\displaystyle{z_{n+1,h}}}{\displaystyle{h}} \\
\ \\ = \frac{\displaystyle{z_{n,h}}}{\displaystyle{h}}
\end{array}
\ \ & \text{in} \ \Omega \\
\ & \ \\
z_{n+1,h}(x) = 0 \ \ & \text{on} \ \partial \Omega \;. \\
\end{array}
\right\}
\label{eq:BasicApprox}
\end{equation}
$z_{n,h}(x)$ is an approximation to $u(x, \; nh),$ and we define the linear
interpolation, $U_h(x,t),$ by
\begin{equation}
U_h(x, \; t) := \frac{[(m+1)h - t]}{h} z_{m,h}(x) + \frac{[t - mh]}{h} z_{m+1,h}(x) \;,
\label{eq:RotheDef}
\end{equation}
where $m \in \N$ is chosen so that $t \in [mh, \; (m+1)h).$
When $U(x, \; t) := \lim_{h \downarrow 0} U_h(x,t)$ exists, we call it the
\textit{Rothe limit}.
 
We make the following continuity assumption, which corresponds to the structure
condition (SC) of \cite{CCKS}:  We assume that for every $R > 0,$ there is a nondecreasing
continuous function $\sigma_R$ such that $\omega_R(0) = 0,$ and if $X$ and $Y$ are symmetric
matrices and $|r|,|s| \leq R,$ then
\begin{equation}
\begin{array}{l}
\hspace{-1in} \ \mathcal{P}^{-}(M - N) - \gamma|P - Q| - \omega_R((s - r)^+) \\
   \leq F(M, \; P, \; r, \; x, \; t) - F(N, \; Q, \; s, \; x, \; t) \\
   \leq \mathcal{P}^{+}(M - N) + \gamma|P - Q| + \omega_R((r - s)^+) \;. \\
\end{array}
\label{eq:ContAss}
\end{equation}
(We will define the Pucci operators $\mathcal{P}^+$ and $\mathcal{P}^-$ in a moment.)
Lastly we assume that for any compact subset of $\R^{n \times n} \times \R^n \times \R$
we have
\begin{equation}
||F(M, \; P, \; r, \; \cdot, \; \cdot)||_{C^{0}(\closure{D})} \leq C < \infty \;.
%||F(0, \; 0, \; 0, \; x, \; t)||_{L^{\infty}(\Omega)} \leq C < \infty \;.
\label{eq:ZeroedStuff}
\end{equation}
 
\newsec{Elliptic Preliminaries}{EllPre}
To study the elliptic problem in each time we introduce the Pucci extremal
operators which we give by
\begin{equation}
\begin{array}{l}
\mathcal{P}^{+}(M) =
    - \lambda \text{Tr}(M^{+}) + \Lambda \text{Tr}(M^{-}) \ \ \ \ \ \ \text{and} \\
\mathcal{P}^{-}(M) =
    - \Lambda \text{Tr}(M^{+}) + \lambda \text{Tr}(M^{-}) \\
%Norm Pucci Defs
%\mathcal{P}^{+}(M) = - \Lambda \; ||M^{-}|| - \lambda \; ||M^{+}|| \ \ \ \ \ \ \text{and} \\
%\mathcal{P}^{-}(M) = - \lambda \; ||M^{-}|| - \Lambda \; ||M^{+}|| \;. \\
\end{array}
\label{eq:PucciXOps}
\end{equation}
where $\text{Tr}(M)$ is the trace of $M.$

\begin{lemma}[First Bound]  \label{FirstBound}
There exists a constant $C$ which is independent of $h$ such that
\begin{equation}
  ||z_{1,h}||_{L^{\infty}} \leq C \;.
\label{eq:LIB}
\end{equation}
\end{lemma}
\begin{pf}
By symmetry it suffices to show that there exists a constant $C$ such that
$z_{1,h}(x) \leq C.$  $z_{1,h}$ satisfies
\begin{alignat*}{1}
0 &= F(D^2 z_{1,h}, \; D z_{1,h}, \; z_{1,h}, \; x, \; h) + \frac{z_{1,h}}{h} \\
  &\geq F(D^2 z_{1,h}, \; D z_{1,h}, \; z_{1,h}, \; x, \; h) \ \ \
   \text{on} \ \ \ \{z_{1,h} > 0 \} \;. \\
\end{alignat*}
and therefore by Lemma 2.11 of \cite{CCKS} $z_{1,h}$ is also
a solution of
\begin{equation}
\mathcal{P}^{-}(D^2 z_{1,h}) - \gamma |D z_{1,h}|
  + F(0, \; 0, \; z_{1,h}, \; x, \; h) \leq 0
\label{eq:CCKSvers}
\end{equation}
on the set where it is positive.  By Proposition 2.12 of \cite{CCKS} we now have
\begin{equation}
z_{1,h} \leq
C||(-F(0, \; 0, \; z_{1,h}, \; x, \; h))^{+}||_{L^{n}(\Gamma^{+}(z_{1,h}))} \;,
\label{eq:fromABP}
\end{equation}
where $\Gamma^{+}$ denotes the upper contact set.
Now by using our structure conditions, we see that
$$-\sigma_R((-z_{1,h})^{+}) \leq F(0, \; 0, \; z_{1,h}, \; x, \; h) -
F(0, \; 0, \; 0, \; x, \; h) \leq \sigma_R((z_{1,h})^{+}) \;.$$
Let $f(x) := F(0, \; 0, \; 0, \; x, \; h).$  Then we see that
\begin{equation}
\left[-F(0, \; 0, \; z_{1,h}, \; x, \; h)\right]^{+}\chisub{ \{ z_{1,h} > 0 \} }
\leq |f(x)| \;.
\label{eq:NextLast}
\end{equation}
By combining Equations\refeqn{fromABP}\!\!,\refeqn{NextLast}\!\!, and\refeqn{ZeroedStuff}we
are done.
\newline Q.E.D. \newline
\end{pf}
\begin{theorem}[First Lipschitz Bound] \label{FLB}
There exists a constant $C$ which is independent of $h$ such that
\begin{equation}
  ||z_{1,h}||_{L^{\infty}} \leq Ch \;.
\label{eq:FLIB}
\end{equation}
\end{theorem}
\begin{pf}
Let the absolute maximum of $z_{1,h}$ be attained at $x_0,$ and let
$\gamma := z_{1,h}(x_0).$  We have
$$\gamma = z_{1,h}(x_0) = -h F(D^2 z_{1,h}(x_0), \; 0, \; \gamma, \; x_0, \; h) \;$$
in the viscosity sense.  Since the plane $\Pi(x) \equiv \gamma$ touches $z_{1,h}$
from above at $x_0$ we know that
\begin{equation}
\gamma \leq -h F(0, \; 0, \; \gamma, \; x_0, \; h)
\label{eq:fa}
\end{equation}
which we know is bounded from above by a constant times $h$ by using the previous
lemma along with Equations\refeqn{ContAss}and\refeqn{ZeroedStuff}\!\!.
\newline Q.E.D. \newline
\end{pf}
 
%Now in general we absolutely cannot expect to have a Lipschitz bound in time.
%We expect $C^{1,\alpha}$ behavior in space but only $(1 + \alpha)/2$-\Holder
%continuity in time.  We still need some sort of estimate of a maximum of the
%difference quotients, but problems arise when we look at their maxima or
%minima, since (as viscosity solutions), they need only be continuous.  For
%this reason, we need to look at mollifications of the $z_n$ and we need to
%find equations satisfied by them.
 
%$$F(z_{n+1},x,(n+1)h) - \frac{z_{n+1}}{h} = - \frac{z_n}{h}$$
\begin{theorem}[Pucci Inequalities]  \label{PEs}
Let $f(x) := F(0,0,0,x,(n+1)h).$  The following inequalities are satisfied
in the viscosity sense.
\begin{equation}
  \begin{array}{l}
     \mathcal{P}^{-}(D^2z_{n+1}) - \gamma||Dz_{n+1}||
       - \sigma_R((-z_{n+1})^{+}) + \displaystyle{\frac{z_{n+1}}{h}} + f
     \leq \ \ \displaystyle{\frac{z_n}{h}} \\
     \ \\
     \leq \mathcal{P}^{+}(D^2z_{n+1}) + \gamma||Dz_{n+1}||
       + \sigma_R(z_{n+1}^{+}) + \displaystyle{\frac{z_{n+1}}{h}} + f \\
  \end{array}
\label{eq:StrIneq}
\end{equation}
\end{theorem}
\begin{pf}
Let $\varphi$ touch $z_{n+1}$ from above at $x_0.$  Since
$$F(D^2 z_{n+1}, D z_{n+1}, z_{n+1},x,(n+1)h) + \frac{z_{n+1}}{h} = \frac{z_n}{h}$$
in the viscosity sense, by using Equation\refeqn{ContAss}we have
\begin{alignat*}{1}
\displaystyle{\frac{z_n}{h}}
   &\geq F(D^2 \varphi, D \varphi, \varphi,x,(n+1)h) + \frac{\varphi}{h} \\
   &\geq \mathcal{P}^{-}(D^2\varphi) - \gamma||D\varphi||
           - \sigma_R((-\varphi)^{+}) + \displaystyle{\frac{\varphi}{h}} + f \;.
\end{alignat*}
The other inequality is proven in the same fashion by touching $z_{n+1}$ from
below.
\newline Q.E.D. \newline
\end{pf}
 
Now we introduce the sup-convolution $z_{j}^{\epsilon}$ of $z_j$ and the
inf-convolution $z_{j,\epsilon}$ of $z_j$ which are defined by
\begin{equation}
   z_{j}^{\epsilon}(x) := \sup_{y \in \Omega}
        \left( z_j(y) - \frac{|x - y|^2}{2 \epsilon} \right) \ \ \ \
   z_{j,\epsilon}(x) := \inf_{y \in \Omega}
        \left( z_j(y) + \frac{|x - y|^2}{2 \epsilon} \right)
\label{eq:SupConvDef}
\end{equation}
The next two lemmas are taken from \cite{J} and \cite{JLS} respectively:
\begin{theorem}[Basic Sup-Convolution Properties] \label{BSCProps}
Let $\Omega \subset \R^n$ be bounded, let $u \in C(\closure{\Omega}),$ and
let $u^{\epsilon}$ and $u_{\epsilon}$ be its sup and inf-convolution
respectively.  Then
  \begin{itemize}
    \item[(i)] $u^{\epsilon},u_{\epsilon} \in C^{0,1}(\Omega).$
    \item[(ii)] $u^{\epsilon} \downarrow u$ and
                $u_{\epsilon} \uparrow u$ as $\epsilon \downarrow 0$
       uniformly on $\closure{\Omega}.$
    \item[(iii)] For every $\epsilon > 0,$ there are measurable functions
       $M^{\epsilon},M_{\epsilon}:\Omega \rightarrow \mathcal{S}(n)$ such that
       $$u^{\epsilon}(y) = u^{\epsilon}(x) + \! < \!Du^{\epsilon}(x),\; y - x \! > \!
         + \frac{1}{2} \! < \!M^{\epsilon}(x)(y - x), \; y - x \! > \! + o(|x - y|^2)$$
       for $a.e. \ x \in \Omega$ (and similarly with $u_{\epsilon}$ and $M_{\epsilon}$)
       where $\mathcal{S}(n)$ is the set of $n \times n$ symmetric matrices.
    \item[(iv)] $M^{\epsilon}(x) \geq -(1/\epsilon)I$ and
                $M_{\epsilon}(x) \leq (1/\epsilon)I$ for $\ a.e \ x \in \Omega.$
    \item[(v)] If $u^{\eta,\epsilon}$ and $u^{\eta}_{\epsilon}$ are standard
       mollifications of $u^{\epsilon}$ and $u_{\epsilon}$ respectively,
       then $D^2 u^{\eta,\epsilon} \geq -(1/\epsilon)I$ and $D^2 u^{\eta,\epsilon}
       \rightarrow M^{\epsilon}(x) \ \ a.e. \ x \in \Omega$ as $\eta \rightarrow 0,$
       and similarly with $D^2 u^{\eta}_{\epsilon}.$
  \end{itemize}
\end{theorem}
\begin{theorem}[Sup-Convolutions are Subsolutions]  \label{SCAsub}
If $u$ is a bounded viscosity subsolution of
$$F(D^2 u, \; Du, \; u, \; x) = f(x) \ \ \text{in} \ B_1$$
and $f, F$ are continuous, then
$$F(M^{\epsilon}(x), \; Du^{\epsilon}(x), \; u(x^{\epsilon}), \; x^{\epsilon})
   \leq f(x^{\epsilon})
   \ \ a.e. \ \ \text{in} \ B_{1 - 2(\epsilon||u||_{\infty})^{1/2}}$$
where $x^{\epsilon} \in B_1$ is any point which satisfies
$$u^{\epsilon}(x) = u(x^{\epsilon}) - \frac{|x^{\epsilon} - x|^2}{2 \epsilon} \;.$$
\end{theorem}
\noindent
Obviously inf-convolutions are supersolutions.
\begin{remark}  \label{OurNot}
What we have denoted by ``$x^{\epsilon},$'' is typically denoted by ``$x^{*},$''
but we need to emphasize the dependence on $\epsilon$ to avoid confusion later.
(For inf-convolutions we will use ``$x_{\epsilon}.$'')
\end{remark}
\begin{remark}  \label{xmxsEst}
Note that
$$\frac{|x^{\epsilon} - x|^2}{2 \epsilon} = u^{\epsilon}(x^{\epsilon}) - u^{\epsilon}(x)
                                   \leq 2||u||_{\infty} \;,$$
so that $|x^{\epsilon} - x| \leq 2(\epsilon ||u||_{\infty})^{1/2}.$  In particular,
if $x \in B_{1 - 2(\epsilon||u||_{\infty})^{1/2}}$ then $x^{\epsilon} \in B_1,$
explaining the appearance of this set above.
\end{remark}
 
By using Equation\refeqn{BasicApprox}and the previous theorem, we get
\begin{equation}
\begin{array}{r}
F(M_{n+1,h}^{\epsilon}(x), \; Dz_{n+1,h}^{\epsilon}(x),
                           \; z_{n+1,h}^{\epsilon}(x^{\epsilon}),
                           \; x^{\epsilon}, \; (n+1)h) +
        \displaystyle{\frac{z_{n+1,h}^{\epsilon}(x^{\epsilon})}{h}} \\
   \ \\
   \leq \displaystyle{\frac{z_{n,h}(x^{\epsilon})}{h}} \\
\end{array}
\label{eq:ViscApprox}
\end{equation}
for $a.e. \; x.$
Similarly we have
\begin{equation}
\begin{array}{r}
F(M_{n,h,\epsilon}(x), \; Dz_{n,h,\epsilon}(x),
                           \; z_{n,h,\epsilon}(x_{\epsilon}),
                           \; x_{\epsilon}, \; nh) +
        \displaystyle{\frac{z_{n,h,\epsilon}(x_{\epsilon})}{h}} \\
   \ \\
   \geq \displaystyle{\frac{z_{n-1,h}(x_{\epsilon})}{h}} \\
\end{array}
\label{eq:ViscApprox2}
\end{equation}
Henceforth we suppress the $h$ in the subscripts.  By taking a difference we get:
\begin{equation}
\begin{array}{r}
F(M_{n+1}^{\epsilon}(x), Dz_{n+1}^{\epsilon}(x),
                            z_{n+1}^{\epsilon}(x^{\epsilon}),
                            x^{\epsilon}, (n+1)h) \\
   \ \\
- \ F(M_{n,\epsilon}(x), Dz_{n,\epsilon}(x),
                           z_{n,\epsilon}(x_{\epsilon}),
                           x_{\epsilon}, nh) \ +
   \displaystyle{\frac{z_{n+1}^{\epsilon}(x^{\epsilon}) -
                       z_{n,\epsilon}(x_{\epsilon})}{h}} \\
   \ \\
   \leq \displaystyle{\frac{z_n(x^{\epsilon}) - z_{n-1}(x_{\epsilon})}{h}} \\
\end{array}
\label{eq:ViscDiff}
\end{equation}
Now we express the difference as a telescoping sum.
\begin{equation}
\begin{array}{r}
F(M_{n+1}^{\epsilon}(x), Dz_{n+1}^{\epsilon}(x),
                            z_{n+1}^{\epsilon}(x^{\epsilon}),
                            x^{\epsilon}, (n+1)h) \\
- \ F(M_{n,\epsilon}(x), Dz_{n,\epsilon}(x),
                           z_{n,\epsilon}(x_{\epsilon}),
                           x_{\epsilon}, nh) \\
   \ \\
= F(M_{n+1}^{\epsilon}(x), Dz_{n+1}^{\epsilon}(x),
                            z_{n+1}^{\epsilon}(x^{\epsilon}),
                            x^{\epsilon}, (n+1)h) \\
- \ F(M_{n,\epsilon}(x), Dz_{n+1}^{\epsilon}(x),
                            z_{n+1}^{\epsilon}(x^{\epsilon}),
                            x^{\epsilon}, (n+1)h) \\
   \ \\
+ \ F(M_{n,\epsilon}(x), Dz_{n+1}^{\epsilon}(x),
                            z_{n+1}^{\epsilon}(x^{\epsilon}),
                            x^{\epsilon}, (n+1)h) \\
- \ F(M_{n,\epsilon}(x), Dz_{n,\epsilon}(x),
                            z_{n+1}^{\epsilon}(x^{\epsilon}),
                            x^{\epsilon}, (n+1)h) \\
   \ \\
+ \ F(M_{n,\epsilon}(x), Dz_{n,\epsilon}(x),
                            z_{n+1}^{\epsilon}(x^{\epsilon}),
                            x^{\epsilon}, (n+1)h) \\
- \ F(M_{n,\epsilon}(x), Dz_{n,\epsilon}(x),
                            z_{n,\epsilon}(x_{\epsilon}),
                            x^{\epsilon}, (n+1)h) \\
   \ \\
+ \ F(M_{n,\epsilon}(x), Dz_{n,\epsilon}(x),
                            z_{n,\epsilon}(x_{\epsilon}),
                            x^{\epsilon}, (n+1)h) \\
%- \ F(D^2 z_{n,\epsilon}(x), Dz_{n,\epsilon}(x),
%                            z_{n,\epsilon}(x_{\epsilon}),
%                            x_{\epsilon}, (n+1)h) \\
%   \ \\
%+ \ F(D^2 z_{n,\epsilon}(x), Dz_{n,\epsilon}(x),
%                            z_{n,\epsilon}(x_{\epsilon}),
%                            x_{\epsilon}, (n+1)h) \\
- \ F(M_{n,\epsilon}(x), Dz_{n,\epsilon}(x),
                            z_{n,\epsilon}(x_{\epsilon}),
                            x_{\epsilon}, nh) \\
\end{array}
\label{eq:Telescope}
\end{equation}
For the time being we will assume that $F$ is differentiable in the
variables we need
%so that $M_{n,\epsilon}$ is simply $D^2 z_{n,\epsilon},$
%and $M_{n+1}^{\epsilon} = D^2 z_{n+1}^{\epsilon}
in order to define
\begin{equation}
\begin{array}{rl}
a_{ij}(x) := \int_{0}^{1} & \!\!\!
             F_{X_{ij}}(tM_{n+1}^{\epsilon}(x) + (1 - t)M_{n,\epsilon}(x), \\
             & \ \ \ Dz_{n+1}^{\epsilon}(x), z_{n+1}^{\epsilon}(x^{\epsilon}),
             x^{\epsilon}, (n+1)h) \; dt \\
\ \\
b_{i}(x) := \int_{0}^{1} & \!\!\!
             F_{P_i}(M_{n,\epsilon}(x), tDz_{n+1}^{\epsilon}(x) +
             (1 - t)D z_{n,\epsilon}(x), \\
             & \ \ \ \ \ \ z_{n+1}^{\epsilon}(x^{\epsilon}),
             x^{\epsilon}, (n+1)h) \; dt \\
\ \\
\tilde{c}(x) := \int_{0}^{1} & \!\!\!
             F_{r}(M_{n,\epsilon}(x), D z_{n,\epsilon}(x),
             tz_{n+1}^{\epsilon}(x^{\epsilon}) + \\
             & \ \ \ (1 - t) z_{n,\epsilon}(x_{\epsilon}),
             x^{\epsilon}, (n+1)h) \; dt \\
\end{array}
\label{eq:CoeffDefs}
\end{equation}
Now we use the fundamental theorem on each pair of terms to get the following
equalities:
\begin{equation}
\left.
\begin{array}{l}
F(M_{n+1}^{\epsilon}(x), Dz_{n+1}^{\epsilon}(x),
                            z_{n+1}^{\epsilon}(x^{\epsilon}),
                            x^{\epsilon}, (n+1)h) \\
- \ F(M_{n,\epsilon}(x), Dz_{n+1}^{\epsilon}(x),
                            z_{n+1}^{\epsilon}(x^{\epsilon}),
                            x^{\epsilon}, (n+1)h) \\
\ \\
= a_{ij}(x) \left(M_{n+1,(ij)}^{\epsilon}(x) - M_{n,\epsilon,(ij)}(x)\right) \\
\end{array}
\right\}
\label{eq:2derEst}
\end{equation}
\begin{equation}
\left.
\begin{array}{l}
F(M_{n,\epsilon}(x), Dz_{n+1}^{\epsilon}(x),
                            z_{n+1}^{\epsilon}(x^{\epsilon}),
                            x^{\epsilon}, (n+1)h) \\
- \ F(M_{n,\epsilon}(x), Dz_{n,\epsilon}(x),
                            z_{n+1}^{\epsilon}(x^{\epsilon}),
                            x^{\epsilon}, (n+1)h) \\
\ \\
= b_{i}(x) D_{i}\left(z_{n+1}^{\epsilon}(x) - z_{n,\epsilon}(x)\right) \\
\end{array}
\right\}
\label{eq:1derEst}
\end{equation}
\begin{equation}
\left.
\begin{array}{l}
F(M_{n,\epsilon}(x), Dz_{n,\epsilon}(x),
                            z_{n+1}^{\epsilon}(x^{\epsilon}),
                            x^{\epsilon}, (n+1)h) \\
- \ F(M_{n,\epsilon}(x), Dz_{n,\epsilon}(x),
                            z_{n,\epsilon}(x_{\epsilon}),
                            x^{\epsilon}, (n+1)h) \\
\ \\
= \tilde{c}(x) \left(z_{n+1}^{\epsilon}(x^{\epsilon}) -
                     z_{n,\epsilon}(x_{\epsilon})\right) \;. \\
\end{array}
\right\}
\label{eq:ConstEst}
\end{equation}
By the continuity assumption\refeqn{ZeroedStuff}we also know
\begin{equation}
\begin{array}{l}
\left| F(M_{n,\epsilon}(x), Dz_{n,\epsilon}(x),
                            z_{n,\epsilon}(x_{\epsilon}),
                            x^{\epsilon}, (n+1)h) \right. \\
\left. - F(M_{n,\epsilon}(x), Dz_{n,\epsilon}(x),
                            z_{n,\epsilon}(x_{\epsilon}),
                            x_{\epsilon}, nh) \right| \\
\ \\
\leq \sigma_1(\epsilon) + \sigma_2(h) \\
\end{array}
\end{equation}
where the $\sigma_i$ are moduli of continuity.  By defining
$c(x) := \tilde{c}(x) + (1/h),$ by defining
$w_{n}^{\epsilon}(x) := z_{n+1}^{\epsilon}(x) - z_{n,\epsilon}(x),$
by defining $D_{ij}w_{n}^{\epsilon} := M_{n+1}^{\epsilon} - M_{n,\epsilon}(x),$
and by using the equations beginning with Equation\refeqn{ViscDiff}we
can conclude
\begin{equation}
\begin{array}{l}
a_{ij}(x) D_{ij}w_{n}^{\epsilon}(x) + b_{i}(x) D_{i}w_{n}^{\epsilon}(x) +
c(x) \left[ w_{n}^{\epsilon}(x^{\epsilon}) + z_{n,\epsilon}(x^{\epsilon}) -
z_{n,\epsilon}(x_{\epsilon}) \right] \\
\ \\
\displaystyle{\leq \frac{z_{n}(x^{\epsilon}) - z_{n - 1}(x_{\epsilon})}{h} +
              \sigma_1(\epsilon) + \sigma_2(h)}
\ \ \ \ \ \text{for} \ a.e. \ x \\
\end{array}
\label{eq:HolyGrail}
\end{equation}
We define a semiconvex function to be a function which becomes convex if
$\gamma ||x||^2$ is added to the function, and $\gamma$ is sufficiently large.
In particular, by Theorem\refthm{BSCProps}\!(iv) we know that sup-convolutions
are semiconvex, as is the function $w_{n}^{\epsilon}$ in Equation\refeqn{HolyGrail}\!\!.
Now we need the following results (see p. 56 Theorem A.2 and p. 58 Lemma A.3
of \cite{CIL}).
\begin{theorem}[Aleksandrov's Theorem]  \label{Aleks}
Semiconvex functions are twice differentiable a.e.
\end{theorem}
\begin{lemma}[Jensen's Lemma]   \label{thmJL}
Let $\varphi: \R^{n} \rightarrow \R$ be semiconvex, and let $\hat{x}$ be a strict
local maximum of $\varphi.$  For $p \in \R^{n},$ set
$\varphi_p(x) = \varphi(x) + <p,x>.$  Then for any $r, \delta > 0,$  the set
$$K := \{ x \in B_r(\hat{x}) \; : \; \text{there exists} \ p \in B_{\delta} \
\text{such that} \ \varphi_p \ \text{has a local max at} \ x \}$$
has positive measure.
\end{lemma}
We also state for future use the following estimate:
\begin{lemma} [Discrete Gronwall Inequality]   \label{DGI}
Assuming $\{v_i\}, \ \{B_i\},$ and $\{D_i\}$ are sequences of nonnegative numbers
which satisfy $v_{i + 1} \leq B_i v_i + D_i,$ we have
\begin{equation}
v_n \leq v_0 \prod_{i=0}^{n - 1} B_i + \sum_{i = 0}^{n - 1}
        \left[ D_i \prod_{j=i + 1}^{n - 1} B_j \right] \;.
\label{eq:DiscGI}
\end{equation}
\end{lemma}
With these results we can turn to the estimate of the maximum of
$w_{n}(x) := z_{n+1}(x) - z_{n}(x).$

\newsec{Rothe Limits Are Lipschitz in Time}{RotheTimeLip}
\begin{theorem}[Inductive Estimation]   \label{IndEst}
If $F$ satisfies\refeqn{ContAss}\!\!, then we have
\begin{equation}
||w_{n}||_{L^{\infty}(\Omega)} \leq
||w_{n-1}||_{L^{\infty}(\Omega)} + h \sigma_2(h) \;.
\label{eq:HappinessAndJoy}
\end{equation}
Furthermore, if we assume that $\sigma_2(h) \leq Ch,$
then for any fixed $T > 0,$ we have the Lipschitz estimate:
\begin{equation}
   \left| \left| \frac{w_{[T/h]}}{h} \right| \right|_{L^{\infty}(\Omega)} \leq C \;,
\label{eq:wnLip}
\end{equation}
where $C$ is independent of $h,$ and $[T/h]$ is the greatest integer less than or equal
to $T/h.$
\end{theorem}
\begin{pf}
We will first prove Equation\refeqn{HappinessAndJoy}\!.
To start we will need a few assumptions:  First that $F$ is differentiable so that
Equation\refeqn{HolyGrail}is valid, second that $w_n$ has a strict maximum at $\hat{x},$
and finally that our approximating functions, $w_n^{\epsilon}$ are twice differentiable
at their corresponding maxima which we call $\hat{x}_{\epsilon}.$  We will show how to
avoid these assumptions after proving the estimate in this case.
 
Because $\hat{x}$ is a maximum, and $w_n^{\epsilon}$ is twice differentiable at
$\hat{x}^{\epsilon},$ Equation\refeqn{HolyGrail}becomes:
\begin{equation}
\begin{array}{l}
c(\hat{x}) \left[ w_{n}^{\epsilon}(\hat{x}^{\epsilon}) +
      z_{n,\epsilon}(\hat{x}^{\epsilon}) -
      z_{n,\epsilon}(\hat{x}_{\epsilon}) \right] \\
\ \\
\displaystyle{\leq \frac{z_{n}(\hat{x}^{\epsilon}) - z_{n - 1}(\hat{x}_{\epsilon})}{h} +
              \sigma_1(\epsilon) + \sigma_2(h)} \\
\end{array}
\label{eq:HG2}
\end{equation}
Now by using Remark\refthm{xmxsEst}we see that
\begin{equation}
|\hat{x}^{\epsilon} - \hat{x}_{\epsilon}| \leq
      2(\epsilon ||z_{n+1}^{\epsilon}||_{\infty})^{1/2}
    + 2(\epsilon ||z_{n,\epsilon}||_{\infty})^{1/2} \;,
\label{eq:DistXs}
\end{equation}
and then by sending $\epsilon$ to zero and using the \Holder regularity of the $z_n$ in
space we can conclude:
\begin{equation}
\tilde{c}(\hat{x})w_{n}(\hat{x}) + \frac{w_{n}(\hat{x})}{h}
          = c(\hat{x}) w_{n}(\hat{x})
       \leq \frac{w_{n-1}(\hat{x})}{h} + \sigma_2(h) \;.
\label{eq:HG3}
\end{equation}
By using Equation\refeqn{ContAss}we see that $F_r \geq 0$ and therefore
$\tilde{c} \geq 0.$  So, since $w_{n}(\hat{x}) \geq 0,$ we have
\begin{equation}
w_{n}(\hat{x}) \leq ||w_{n-1}||_{L^{\infty}(\Omega)} + h\sigma_2(h) \;.
\label{eq:QED}
\end{equation}
 
Obviously we can argue similarly for a negative minimum, but it remains to show how
to eliminate our regularity assumptions.  First of all we need stability of our
approximations as we take smooth approximators $F_m,$ to our original $F.$  The
desired stability can be found in section 6 of \cite{CIL}.  (See Remark 6.3 in
\cite{CIL} in particular.)  In terms of applying the estimates we conclude for our
$F_m$'s for our original $F$ there is no problem, because no matter how badly
behaved the $a_{ij}, \; b_i,$ and $\tilde{c}$ are in a given approximation to $F,$
they get thrown away at the maximum of $w_{n}$ by considerations which depend only
on the signs of the terms involved.  Stated differently, even if the derivatives of
the approximations are diverging, they are necessarily diverging ``in the right
direction.''  (Equation\refeqn{HG2}is obtained from Equation\refeqn{HolyGrail}by
using the nonpositivity of the second derivative and the vanishing of the gradient
of $w_n$ at its maximum.)
 
Next, if $w_{n}$ has a maximum at $\hat{x}$ but it is not the unique
maximum, then we look at the equation satisfied by
$\tilde{w}_{n}(x) := w_{n}(x) - \hat{\epsilon}|x - \hat{x}|^2 \;,$
where $\hat{\epsilon} > 0$ is arbitrarily small.  (Here we need the continuity
assumptions on our function $F$ given in\refeqn{ContAss}\!\!.)
Finally, to ensure twice
differentiability of the $w_{n}^{\epsilon}$ at $\hat{x}^{\epsilon}$ we appeal to
Jensen's Lemma\refthm{thmJL}and Aleksandrov's Theorem\refthm{Aleks}\!\!.

Now we need to prove Equation\refeqn{wnLip}\!\!.
We apply Lemma\refthm{DGI}to Equation\refeqn{HappinessAndJoy}to get
\begin{equation}
||w_{n}||_{L^{\infty}} \leq ||w_{0}||_{L^{\infty}} + nh \sigma_2(h) \;.
\label{eq:LipEst1}
\end{equation}
with $n = [T/h]$ and by using Theorem\refthm{FLB}we are done.
\newline Q.E.D. \newline
\end{pf}
\begin{theorem}[Existence of the Rothe Limit] \label{RotheExists}
If $\sigma_2(h) \leq Ch$ the Rothe limit exists and is locally Lipschitz
in time.
\end{theorem}
\begin{pf}
The last theorem shows that the sequence of approximators are uniformly Lipschitz
in time.  In order to show that the approximators are equicontinuous in space, we
observe that our uniform Lipschitz estimate implies
$F(D^2 z_{n+1}, D z_{n+1}, z_{n+1}, x, (n+1)h)$ is uniformly bounded in $L^{\infty}$
in the viscosity sense.  Then we invoke the spatial \Holder estimates known for fully
nonlinear elliptic equations (see Theorem 5.1 Part III of \cite{T}).  Finally, we can
invoke Arzela-Ascoli to guarantee the existence of the limit.
\newline Q.E.D. \newline
\end{pf}
 
\newsec{Rothe Limits Are Viscosity Solutions}{RotheVisc}
\begin{theorem}[Rothe Limits Are Viscosity Solutions]   \label{RotheisVisc}
The Rothe limit is a viscosity solution.
\end{theorem}
\begin{pf}
By symmetry, it will suffice to prove that the Rothe limit is a viscosity
subsolution.  Suppose not.  Then there exists an $\epsilon > 0$ and a supersolution
$\varphi \in C^{2,1}(\closure{D})$ which touches the limit $U$ from above at a
point $(x_0, t_0)$ and which satisfies:
     \begin{equation}
     \varphi_t + F(D^2 \varphi, \; D \varphi, \; \varphi, \; x, \; t) \geq \epsilon > 0 \;,
     \label{eq:vpissuper2}
     \end{equation}
pointwise in a neighborhood of $(x_0, t_0).$
By adding a very small multiple of $|t - t_0|^2 + |x - x_0|^2$ to $\varphi,$ by using
the continuity of $F,$ and by allowing a slightly smaller (but still positive)
$\epsilon$ in Equation\refeqn{vpissuper2}\!\!, we can assume without loss of
generality that at least for a small neighborhood of $(x_0,t_0)$ which we will
call $N$ the only contact between $\varphi$ and $U$ is at $(x_0, t_0).$  Because
of the uniform continuity of $\varphi_t$ there exists a neighborhood $\tilde{N}$
of $(x_0, t_0)$ such that for a sufficiently small $\tilde{h} > 0$ and for all
$(x, t) \in \tilde{N}$ we have the estimate
\begin{equation}
   \left| \frac{\varphi(x,t) - \varphi(x,t-h)}{h} - \varphi_t(x,t) \right| \leq
          \frac{\epsilon}{10}
\label{eq:unifderest}
\end{equation}
for all $h \leq \tilde{h}.$  At this point let $\mathcal{N}$ denote a
compact set of the form $\qclosure{B_{r}(x_0)} \times [t_0 - \gamma, t_0 + \gamma]$ which is
contained in $N \cap \tilde{N}.$  Let
\begin{equation}
\mathcal{S} := \mathcal{N} \setminus \{ B_{r/2}(x_0) \times
   (t_0 - \gamma/2, t_0 + \gamma/2) \}.
\label{eq:Sdef}
\end{equation}
Because $\varphi > U$ on $\mathcal{S},$ there exists a $\tilde{\epsilon} > 0$
such that $\varphi \geq U + \tilde{\epsilon}$ on $\mathcal{S}.$
By the uniform convergence of $U_h$ to $U,$ for any
$\delta \in (0,\tilde{\epsilon}/3)$ we can be sure that by shrinking $\tilde{h}$
if necessary we have
\begin{equation}
   \varphi > U_h + \tilde{\epsilon} - \delta \ \ \text{on} \ \mathcal{S}
   \ \ \text{while} \ \ | \; \varphi(x_0,t_0) - U_h(x_0,t_0) | < \delta/4 \;.
\label{eq:stillok}
\end{equation}
for all $h \leq \tilde{h}.$  Choose $\delta$ sufficiently small to guarantee that
\begin{equation}
   \omega_R(\delta) < \epsilon/2 \;,
\label{eq:delisgood}
\end{equation}
where $\omega_R$ is the modulus given in Equation\refeqn{ContAss}and $\epsilon$ is
taken from Equation\refeqn{vpissuper2}\!\!.
 
Fix $h < \min \{\tilde{h}, \gamma/10\},$ and sufficiently small to guarantee that
\begin{equation}
| \varphi(x_0,t) - U_h(x_0,t) | < \delta/2
\label{eq:minonmesh}
\end{equation}
for any $t \in [t_0 - h,t_0].$  (Here we need the second half of
Equation\refeqn{stillok}and uniform continuity in time of
$\varphi(x_0,t) - U_h(x_0,t)$ as $h \downarrow 0.$  The continuity is guaranteed by
the previous theorem.)  If we let $\tilde{\varphi} = \varphi + c,$
for any $c$ with $|c| < \delta,$ then by using
Equations\refeqn{ContAss}and\refeqn{delisgood}we will have
\begin{equation}
\tilde{\varphi}_t + F(D^2 \tilde{\varphi}, \; D \tilde{\varphi}, \; \tilde{\varphi},
   \; x, \; t) \geq \epsilon - \omega_R(\delta) > \epsilon/2 > 0 \;.
\label{eq:vtilstsup}
\end{equation}
Now we choose $c$ so that the minimum of $\tilde{\varphi} - U_h$ on the set
$\mathcal{N} \cap \{ (x,t) : t \in \N h \}$ is equal to zero.  In other words,
$\tilde{\varphi}$ is allowed to be less than $U_h$ between the mesh values of
$t,$ but when considering \textit{only} mesh values of $t,$ we can say that
$\tilde{\varphi}$ touches $U_h$ from above.
 
\
 
\psfig{file=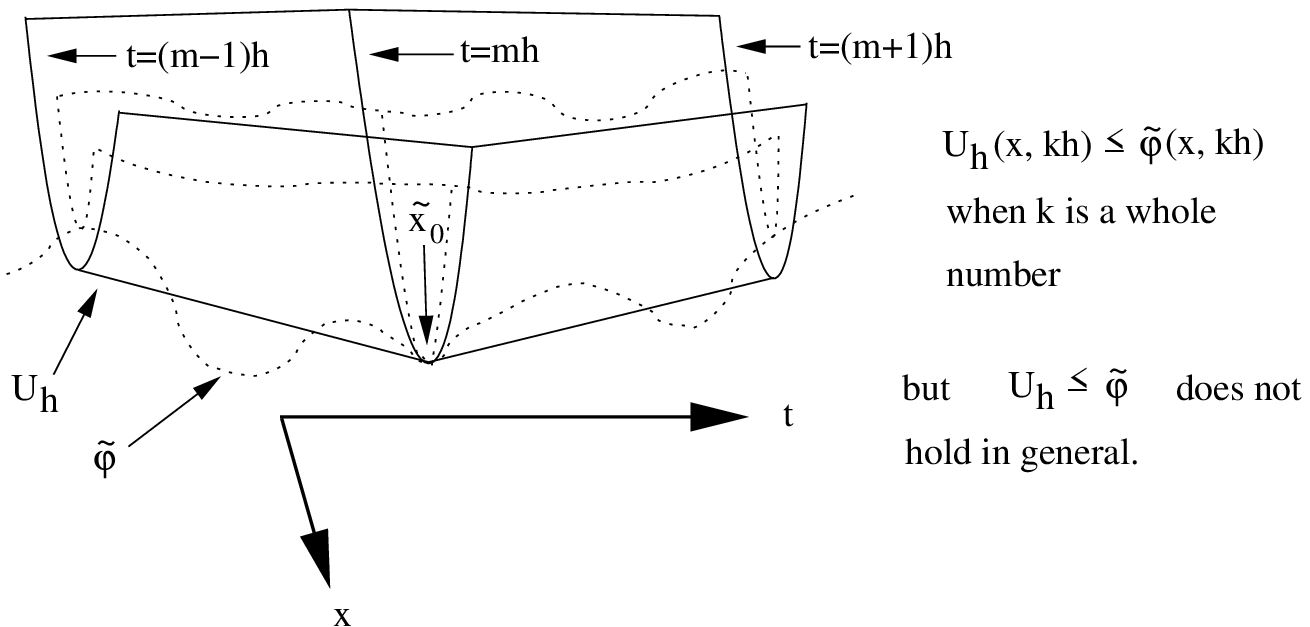}
 
\
 
\noindent
Because of Equation\refeqn{minonmesh}\!,
we can be sure that $|c| \leq \delta/2.$  Let $(\tilde{x_0},mh)$ be a point
where $\tilde{\varphi} = U_h.$  By using Equation\refeqn{stillok}we know that
$(\tilde{x_0},mh) \in B_{r/2}(x_0) \times (t_0 - \gamma/2, t_0 + \gamma/2).$
 
If we return to our $z_n$ we now have the following setting:
\begin{itemize}
   \item[(a)] $\tilde{\varphi}(x,mh)$ touches $z_m(x)$ from above at $\tilde{x_0}.$
   \item[(b)] $\tilde{\varphi}(x,(m-1)h) \geq z_{m-1}(x)$ in $B_{r/2}(x_0).$  (Because
              $h < \gamma/10$ we know that $B_{r/2}(x_0) \times \{ t = (m-1)h \} \subset
              \mathcal{N}.$)
\end{itemize}
Because $(\tilde{x_0},mh) \subset \tilde{N}$ we have (see Equation\refeqn{unifderest}\!)
\begin{equation}
   \left| \frac{\tilde{\varphi}(\tilde{x_0},mh) -
                \tilde{\varphi}(\tilde{x_0},(m-1)h)}{h} -
                \tilde{\varphi}_t(\tilde{x_0},mh) \right| \leq
          \frac{\epsilon}{10}
\label{eq:diffquotgood}
\end{equation}
Recall that $z_m$ is a viscosity solution of
\begin{equation}
   \frac{z_m - z_{m - 1}}{h} + F(D^2 z_m, D z_m, z_m, x, mh) = 0 \;.
\label{eq:zmvisc}
\end{equation}
Since $\tilde{\varphi}(x,mh)$ touches $z_m(x)$ from above at $\tilde{x_0},$ we
conclude that
\begin{equation}
   \frac{\tilde{\varphi} - z_{m - 1}(\tilde{x_0})}{h} +
   F(D^2 \tilde{\varphi}, D \tilde{\varphi}, \tilde{\varphi}, \tilde{x_0}, mh) \leq 0 \;.
\label{eq:varphivisc}
\end{equation}
Now we combine Equations\refeqn{vtilstsup}\!, \refeqn{varphivisc}\!,
and\refeqn{diffquotgood}to get
\begin{alignat*}{1}
     \epsilon/2 \;
  &< \; \tilde{\varphi}_t(\tilde{x_0}) +
        F(D^2 \tilde{\varphi}, D \tilde{\varphi}, \tilde{\varphi}, \tilde{x_0}, mh) \\
  &\leq \; \tilde{\varphi}_t(\tilde{x_0}) +
        \frac{z_{m - 1}(\tilde{x_0}) - \tilde{\varphi}(\tilde{x_0}, mh)}{h} \\
  &\leq \; \tilde{\varphi}_t(\tilde{x_0}) +
        \frac{\tilde{\varphi}(\tilde{x_0}, (m - 1)h) - \tilde{\varphi}(\tilde{x_0}, mh)}{h} \\
  &\leq \epsilon/10
\end{alignat*}
which is a contradiction.
\newline Q.E.D. \newline
\end{pf}
\begin{corollary}[Lipschitz Regularity in Time]  \label{LipTime}
There is a solution to Equation\refeqn{BasicEqn}which is Lipschitz in time.
\end{corollary}
\begin{pf}
Simply combine the two previous results.
\newline Q.E.D. \newline
\end{pf}

\end{document}